\documentclass{amsart}
\usepackage{latexsym,amssymb,amsthm,amsmath}

\usepackage{amssymb}
\usepackage{amsmath}

\theoremstyle{plain}
\newtheorem{theorem}{Theorem}[section]
\newtheorem*{Theorem B}{Theorem B}
\newtheorem*{Theorem A}{Theorem A}
\newtheorem{lemma}{Lemma}[section]

\newtheorem{corollary}{Corollary}[section]

\numberwithin{equation}{section}

\theoremstyle{remark}
\newtheorem{remark}{Remark}[section]
\newtheorem{acknowledgements}{Acknowledgement}
\title[B.-Y. Chen's inequality  for  $CR$-warped products]{B.-Y. Chen's inequality  for  $CR$-warped products in a locally conformal Kaehler space form}

\author{Siraj Uddin, Khushwant Singh, Fatimah Alghamdi}
\address{S. Uddin: Department of Mathematics, Faculty of Science, King Abdulaziz University, 21589 Jeddah, Saudi Arabia}
\email{siraj.ch@gmail.com}
\address{K. Singh:  Indian Institute of Science Education and Research, S.A.S. Nagar Mohali-146306, India}
\email{khushwantchahil@gmail.com}
\author{Cengizhan Murathan}
\address{F. Alghamdi: Department of Mathematics, Faculty of Science, Jeddah University, 21589 Jeddah, Saudi Arabia}
\email{fmalghamdi@uj.edu.sa}
\address{C. Murathan: Department of Mathematics, Uludag University, 16059 Bursa, Turkey}
\email{ gengiz@uludag.edu.tr}


\begin{document}
\begin{abstract} 
In this paper, we obtain a geometric inequality between the length of the second fundamental form and the length of Lee form in terms of the warping function for a CR-warped product submanifold in a locally conformal Kaehler space form. The equality case is also investigated. Furthermore, the inequality is discussed for the important subclass of locally conformal Kaehler manifolds i.e., Vaisman manifold.  
\end{abstract}

\maketitle

\sloppy
\section{Introduction}
It is well-known that the notion of warped products plays some important role in differential geometry as well as physics. R. L. Bishop and B. O'Neill in 1969
introduced the concept of a warped product manifold to provide a
class of complete Riemannian manifolds with everywhere negative
curvature \cite{Bi}. The warped product scheme was later applied to
semi-Riemannian geometry (for instance, see \cite{JK}) and general relativity \cite{JK2}.

Recently, Chen \cite{Chen01} (see also \cite{Chen001}) studied warped product, he considered warped product
$CR$-submanifolds in the form $M=N^T\times{_{f}N^\perp}$ for which he
called $CR$-warped products, where $N^T$ and $N^\perp$ are
holomorphic and totally real submanifolds of a Kaehler manifold
$\tilde{M}$. Later, he extended the study of CR-warped products in real and complex space forms (see, \cite{C6, C7}. Motivated by Chen's papers many authors studied
$CR$-warped product submanifolds in almost complex as well as
contact setting (see \cite{Bo}, \cite{CU} \cite{Ja10}, \cite{Ha03}, \cite{Mun05,Mun07}, \cite{SU,SU1, SC}). For up-to-date surgery on the warped product manifolds and warped product submanifolds we refer Chen's books \cite{book, book17}.

In this paper, we have obtained a
general sharp inequality for the length of second fundamental form
of $CR$-warped product submanifolds in a locally conformal Kaehler
space form (in short $LCK$-space form). Also, the inequality is
discussed for a Vaisman manifold.


\section{Preliminaries}
A locally conformally Kaehler $(LCK)$ manifold $M$ is one which is
covered by a Kaehler manifold $\tilde{M}$ with the deck
transformation group acting conformally on $\tilde{M}$. $LCK$
manifolds have been widely studied in the last 30 years (see
(see, \cite{Mun05,Mun07}, \cite{Or05,Or10}, \cite{Va76,Va80}). They share some properties with Kaehler manifolds  c.f. \cite{Or05}.

A Hermitian manifold $\tilde {M}$ with structure $(J,\tilde{ g})$ is called a locally conformal
Kaehler (an $ LCK $ ) manifold if each point $x \in M$ has an open neighbourhood $U$ with differentiable function $ \rho : U \rightarrow R$ such that $\tilde {g}^{*}= e^{-2\rho}\tilde g\mid U$ is a Kaehlerian   metric  on  $U$, that  is, $\nabla^{*}J = 0$, where $J$ is the almost complex
structure, $\tilde{g}$ is the Hermitian metric, $\nabla^{*}$
is the covariant differentiation with
respect to $\tilde{g}$ and $R$ is a real number space.

Let $\tilde{M}$ be an $LCK$ manifold. Then the vector
field $\lambda$ (the Lee field of $\tilde{M}$) is defined by
$g(X,\lambda) = \alpha(X)$. The best known examples of $LCK$ manifolds are the Hopf manifolds.

\begin{theorem}\label{T1} \cite{Va76} The almost Hermitian manifold $\tilde{M}$ is an $LCK$ manifold if and only if there is a closed $1$-form $ \alpha$ on $\tilde{M}$ which is called Lee form.
\end{theorem} If $\tilde{\nabla}$ denotes the Levi-Civita connection on
$\tilde{M}$ , then we have
\begin{align}\label{2.1}
(\tilde{\nabla}_U J)V = -g(\beta ^{\#}, V )U -g(\alpha^{\#}, V )JU + g(JU, V )\alpha ^{\#} +g(U, V )\beta^{\#} .
\end{align}
for any $ U, V $ $\in$ $\tilde {TM}$, where $ \alpha^{\#}$ is the dual vector field of $ \alpha$ which is called the Lee vector field, $\beta $ is the 1-form defined by
$\beta(U) =\alpha(JU)$ for any $U$ $\in$ $ \tilde {TM}$ and $\beta^{\#}$ is the dual vector field of $\beta$ \cite{Va76}. In
terms of the Lee vector field, above equation can be written as
\begin{align}\label{2.2}
(\tilde{\nabla}_U J)V = [g(\lambda, JV )U-g(\lambda, V )JU + g(JU, V )\lambda + g(U, V )J\lambda].
\end{align}
The most important subclass of $LCK$ manifolds is Vaisaman manifold and it is defined by the
parallelism of the Lee form with respect to the Levi-Civita
connection of $g$. An $LCK$ manifold $(\tilde{M}, J, g)$
is called a Vaisman manifold if $\tilde\nabla \alpha = 0$, where $\tilde\nabla$
is the Levi-Civita connection of $g$ (see, for instance \cite{Va76,Va80}).

An $LCK$-manifold $\tilde{M}$ is called an $LCK$-space form if
it has a constant holomorphic sectional curvature $c$. Then the
Riemannian curvature tensor. $\tilde{R}$ of, an $LCK$-space form
$\tilde{M}(c)$ with constant holomorphic sectional curvature $c$ is
given by is given by Matsumoto \cite{Mat84}
\begin{align}\label{2.3}
\tilde{R}(X,Y,Z,W)&=\frac{c}{4}\{g(X,W)g(Y,Z)-g(X,Z)g(Y,W)+g(JX,W)g(JY,Z)\notag\\
&-g(JX,Z)g(JY,W)-2g(JX,Y)g(JZ,W)\}+\frac{3}{4}\{\ P(X,W)g(Y,Z)\notag\\
&-P(X,Z)g(Y,W)+g(X,W)P(Y,Z)-g(X,Z)P(Y,W)\}\notag\\
&-\frac{1}{4}\{\tilde P(X,W)g(JY,Z)-\tilde P(X,Z)g(JY,W)+g(JX,W)\tilde P(Y,Z)\notag\\
&-g(JX,Z)\tilde P(Y,W)\}+\frac{1}{2}\{\ \tilde P(X,Y)g(JZ,W)+\tilde P(Z,W)g(JX,Y)\}
\end{align}
where $\tilde{R}(X,Y,Z,W)=g(\tilde{R}(X,Y)Z,W)$ and $P,\,\,\tilde P$ respectively defined by 
\begin{align}\label{2.4}
P(Y,X)=-(\tilde{\nabla}_Y \alpha)X-\alpha(Y)\alpha(X) +\frac{1}{2}\| \alpha\|^2g(X, Y),
\end{align}
where $\| \alpha\|^2$ denotes the length of the Lee form $\alpha$ with respect to $g$ and
\begin{align}\label{2.4a}
\tilde P(X,Y)=P(JX,Y)
\end{align}
\begin{remark}\label{R1} $\tilde P(X,Y)$ is skew-symmetric.
\end{remark}
Let $\tilde{M}(J,g,\alpha)$ be a complex $m$-dimensional
$LCK$-manifold and $M$ be a real $n$-dimensional ($n\leq m$) Riemannian manifold
isometrically immersed in $\tilde{M}$. We denote the metric tensor
induced on $M$ by $g$. Let $\nabla$ be the covariant differentiation
with respect to the induced metric on $M$. Then the Gauss and
Weingarten formulas for $M$ are respectively given by
\begin{align}\label{2.5}
\tilde {\nabla}_X Y={\nabla}_X Y+ h(X,Y),\;\;\;\tilde {\nabla}_X\xi=-A_\xi X+ {\nabla^\perp}_X\xi,
\end{align}
for any $X,Y$ tangent to $M$ and $\xi$ normal to $M$, where
$\nabla^\perp$ is the connection on the normal bundle $T^\perp M $,
$h$ is the second fundamental form and $A_N$ is the Weingarten map associated with the vector field $N\in T^\perp M$ as
\begin{align}\label{2.6}
g(A_\xi X,Y)=g(h(X,Y), \xi).
\end{align}
The covariant derivative of the second fundamental form is given by
\begin{align}\label{2.7}
(\tilde{(\nabla}_X h)(Y,Z)=\tilde\nabla_X h(Y,Z)-h(\nabla_X Y,Z)-h(Y,\nabla_X Z),
\end{align}
for all $X,Y, Z \in TM$.

A Riemannian manifold $M$, isometrically immersed in an $LCK-$manifold $\tilde M$ is called a {\it{CR-submanifold}} if
there exist on $M$ a differentiable holomorphic distribution $\mathfrak{D}$
i.e., $J\mathfrak{D}_x=\mathfrak{D}_x$, for any $x\in M$ whose orthogonal complementary distribution $\mathfrak{D}^\perp$ in $TM$ is totally real on $M$ i.e., $J\mathfrak{D}_x^\perp\subset
T_x^\perp M$. For a CR-submanifold $M$ of an $LCK-$manifold $\tilde{M}$, the normal bundle $T^\perp M$ is decomposed as
\begin{align*}T^\perp M=J\mathfrak{D}^\perp\oplus\nu,
\end{align*}
where $\nu$ is the invariant normal subbundle of $T^\perp M$ under $J$. Now, on a CR-submanifold of an $LCK-$manifold $\tilde{M}$, we have the following useful result.
\begin{lemma}\label{L1} Let $M$ be a $CR$-submanifold of an $LCK$-manifold
$\tilde{M}$. Then we have
\begin{enumerate}
\item[{(i)}] $g(\nabla_U Z,X)=g(JA_{JZ}U,X)-g(J\lambda,Z)g(JU,X)-g(U,Z)g(\lambda,X)+g(\lambda,Z)g(U,X),$
\item[{(ii)}] $A_{JZ}W-A_{JW}Z=g(J\lambda,Z)W-g(J\lambda,W)Z,$
\item[{(iii)}] $ A_{J\xi}X+A_{\xi} JX=g(J\lambda,\xi)X-g(\lambda,\xi)JX+g(\lambda,JX)\xi+g(\lambda,X)J\xi,$
\end{enumerate}
for any $X\in\mathfrak{D},\,\,Z,W\in \mathfrak{D}^\perp;\;\;\xi\in\nu$ and $U\in TM.$
\end{lemma}
\begin{proof} The proof is straightforward and obtained by using \eqref{2.2}, \eqref{2.5} and \eqref{2.6}.
\end{proof}
Let us calculate the holomorphic bisectional curvature
$\tilde{H}_B (X, Z)$ for unit vectors $X \in\mathfrak{D}$ and $Z \in \mathfrak{D}^\perp$,
where $\tilde{H}_B(X,Z)$ is defined by
\begin{align*}
\tilde{H}_B(X, Z)=\tilde{R}(X, JX;JZ,Z).
\end{align*}
By the straightforward calculation, we get the following
lemma.

\begin{lemma}\label{L2}
Let $\tilde{M}$ be an $LCK$-space form and let $X\in\mathfrak{D}$ and $Z\in
\mathfrak{D}^\perp$ be unit vector fields. Then the holomorphic bisectional
curvature of the plane $X\wedge Z$ is given by
\begin{align}\label{2.8}
\tilde{H}_B (X,Z)=\frac{c}{2}-\frac{1}{2}\left\{P(X,X)+P(Z,Z)\right\}. 
\end{align}
\end{lemma}

\begin{proof} By definition, we know that
\begin{align*}
\tilde{H}_B(X, Z)=\tilde{R}(X, JX;JZ,Z).
\end{align*}
By using equation \eqref{2.3}, we get
\begin{align*}
\tilde{R}(X,JX,JZ,Z)&= \frac{c}{4}\{g(X,Z)g(JX,JZ)-g(X,JZ)g(JX,Z)\notag\\
&-g(JX,Z)g(X,JZ)+g(JX,JZ)g(X, Z)+2g(JX,JX)g(Z,Z)\}\notag\\
&+\frac{3}{4}\{P(X,Z)g(JX,JZ)-P(X,JZ)g(JX,Z)+g(X,Z)P(JX,JZ)\notag\\
&-g(X,JZ)P(JX,Z)\}+\frac{1}{4}\{\tilde P(X,Z)g(X,JZ)-\tilde P(X,JZ)g(X,Z)\notag\\
&-g(JX,Z)\tilde P(JX,JZ)+g(JX,JZ)\tilde P(JX,Z)\}\notag\\
&-\frac{1}{2}\{\tilde P(X,JX)g(Z,Z)-\tilde P(JZ,Z)g(JX,JX)\},
\end{align*}
for any $X\in\mathfrak{D}$ and $Z\in \mathfrak{D}^\perp$ in the plane $X\wedge Z$, from
above equation it follows that
\begin{align*}\tilde{H}_B (X,Z)=\frac{c}{2}+\frac{1}{2}\left(\tilde P( JZ,Z)- \tilde  P(X,JX)\right).
\end{align*}
By virtue of \eqref{2.4a} and Remark \ref{R1}, we obtain the required result.
\end{proof}

In case of Vaisman manifold, from above lemma we get the
following important result.
\begin{corollary}\label{C1}
Let $\tilde{M}$ be a Vaisman manifold. Then the holomorphic
bisectional curvature of the plane $X\wedge Z$ is given by
\begin{align}\label{2.9}
\tilde{H}_B (X,Z)=\frac{c}{2}+\frac{1}{2}\{g(Z,\lambda)^2+ g(X,\lambda)^2-\|\alpha\|^2\}
\end{align}
where $\lambda$ is the Lee vector field.
\end{corollary}
\begin{proof}
The proof follows from \eqref{2.2} and \eqref{2.8}.
\end{proof}

\section{Warped product CR-submanifolds}

In 1967, Bishop and O'Neill introduced the notion of warped product
manifolds \cite{Bi}. They defined these manifolds as: Let $(N_1, g_1)$ and
$(N_2, g_2)$ be two Riemannian manifolds and $f>0$ a differentiable
function on $N_1$. Consider the product manifold $N_1\times N_2$
with its projections $\pi_1:N_1\times N_2\to N_1$ and
$\pi_2:N_1\times N_2\to N_2$. Then the warped product of $N_1$ and
$N_2$ denoted by $M=N_1\times{_f}N_2$ is a Riemannian manifold
$N_1\times N_2$ equipped with the Riemannian structure such that
$$g(X, Y)=g_1({\pi_1}_\star X, {\pi_1}_\star Y)+(f\circ\pi_1)^2g_2({\pi_2}_\star X, {\pi_2}_\star Y)$$
for each $X, Y\in\Gamma(TM)$ and $\star$ is a symbol for the tangent
map. Thus we have
\begin{align*}
g=g_1+f^2g_2.
\end{align*}
The function $f$ is called the {\it{warping function}} of the warped
product \cite{Bi}. A warped product manifold $N_1\times{_f}N_2$ is said to
be {\it{trivial}} if the warping function $f$ is constant.

We recall the following general result obtained by Bishop and O'Neill \cite{Bi} for warped product manifolds.

\begin{lemma}\label{L3}\cite{Bi} Let $M=N_1\times{_{f} N_2}$ be a warped product manifold with
the warping function $f$, then for any $X, Y\in TN_1$ and $Z, W\in TN_2$, we have
\begin{enumerate}
\item [{(i)}] $\nabla_XY\in T(N_1)$,
\item [{(ii)}] $\nabla_XZ=\nabla_ZX=(X\ln f)Z$,
\item [{(iii)}] $\nabla_ZW=\nabla_Z^{N_2}W-g(Z,W)\nabla\ln f$,
\end{enumerate}
where $\nabla$ and $\nabla^{N_2}$ denote the Levi-Civita
connections on $M$ and $N_2$, respectively and $\nabla\ln f$ is the
gradient of the function $\ln f$.
\end{lemma}

In this section, we study CR-warped product submanifolds of the form $N^T\times_fN^\perp$ in LCK-space forms. For the simplicity, throughout this paper we denote the corresponding tangent spaces of $N^T$ and $N^\perp$ by $\mathfrak{D}$ and $\mathfrak{D}^\perp$, respectively.

\begin{lemma}\label{L4} Let $M=N^T\times_f N^\perp$ be a $CR$-warped product submanifold
in an $LCK$-manifold $\tilde{M}$. Then
\begin{enumerate}
\item [{(i)}] $g(h(X,Y),JZ)=g(J\lambda,Z)g(X,Y)$,
\item [{(ii)}] $g(h(JX,W),JZ)=(X \ln f-g(\lambda,X))g(Z,W)$,
\item [{(iii)}] $g(h(X,W),JZ)=-(JX \ln f+g(J\lambda,X))g(Z,W)$,
\end{enumerate}
for any $X, Y$ tangent to $N^T$ and $Z,W$ tangent to $N^\perp.$
\end{lemma}
\begin{proof} statements (i) is proved in \cite{Bo} (see Proposition 3.1) but For (ii), by Lemma \ref{L3} $(ii)$, we have
\begin{align*}
\nabla_X Z=(Z\ln f)X
\end{align*}
for any $X\in T(N^T)$ and $Z\in T(N^\perp)$. Then from \eqref{2.2}, we get
\begin{align*}g(h(JX,W),JZ)&=g(\lambda,JX)g(W,JZ)-g(\lambda,X)g(JW,JZ)+g(JW,X)g(\lambda,JZ)\notag\\
&+g(X,W)g(J\lambda,JZ)+(X\ln f)g(W,Z).
\end{align*}
Thus, the result follows from the above relation.
\end{proof}
Now we will compute the norm of the $\nu$-component  of $h(X,Z)$ and norm of the $J\mathfrak{D}^\perp$-component of  $h(X,Z)$. 

\begin{lemma}\label{L5} Let $M=N^T\times_f N^\perp$ be a  $CR$-warped product submanifold of an LCK-manifold $\tilde M$. Then, we have
\begin{enumerate}
\item [{(i)}] $\|h_\nu(X,Z)\|^2=g(Jh(X,Z),h(JX,Z)),$
\item [{(ii)}] $\|h_{J\mathfrak{D}^\perp}(X,Z)\|^2=(X \ln f-g(\lambda,X))^{2}\|Z\|^2.$
\end{enumerate}
for any $X, Y\in\mathfrak{D}$ and $Z, W\in\mathfrak{D}^\perp$.
\end{lemma}

\begin{proof}
For the first part of this lemma we know that
\begin{align*}
\|h_\nu(X,Z)\|^2=g(h_\nu(X,Z),h(X,Z))=g(A_{h_\nu (X,Z)} X,Z).
\end{align*}
By using Lemma \ref{L1} $(iii)$ , we get
\begin{align*}\|h_\nu(X,Z)\|^2& =g(Jh_\nu(X,Z),h(JX,Z))+g(J\lambda,Jh_\nu(X,Z))g(X,Z)\\ +
&g(\lambda,Jh_\nu(X,Z))g(JX,Z)-g(\lambda,JX)g(Z,Jh_\nu(X,Z))\\
&+g(\lambda,X)g(Z,h_\nu(X,Z)).
\end{align*}
for any $X\in\mathfrak{D}$ and $Z\in \mathfrak{D}^\perp$. From
above relation it follows that
\begin{align*} 
\|h_\nu(X,Z)\|^2&=g(Jh_\nu(X,Z),h(JX,Z))\\
&=g(Jh(X,Z),h(JX,Z))-g(Jh_{J\mathfrak{D}^\perp}(X,Z)\,h(JX,Z)).
\end{align*}
From the fact that $Jh_{J\mathfrak{D}^\perp}(X,Z)$
belongs to $\mathfrak{D}^\perp$, we obtain
\begin{align*}
\|h_\nu(X,Z)\|^2&=g(Jh(X,Z),h(JX,Z))
\end{align*}
which is (i). For the second part of this lemma we know that
\begin{align*}
\|h_{J\mathfrak{D}^\perp}(X,Z)\|^2=-g((h(JX,Z)\,J^2 h_{J\mathfrak{D}^\perp}(JX,Z)).
\end{align*}
By using Lemma \ref{L4} $(iii)$, we get 
\begin{align*}
\|h_{J\mathfrak{D}^\perp}(X,Z)\|^2&=-(X  lnf-g(\lambda,X))g(Jh_{J\mathfrak{D}^\perp}(JX,Z),Z)\notag\\
&=(X\ln f-g(\lambda,X))g(h_{J\mathfrak{D}^\perp}(JX,Z),JZ)\notag\\
&=(X\ln f-g(\lambda,X))g(h(JX,Z),JZ)\notag\\
&=(X\ln f- g(\lambda,X))^{2}\|Z\|^2,
\end{align*}
which is (ii). Hence, the lemma is proved.
\end{proof}

The following  characterization is proved as a Proposition 3.1 in \cite{Bo} and Theorem 3.5 of \cite{Ja10}, "a $CR$-submanifold
$M$ of an $LCK$-manifold $\tilde{M}$ is a $CR$-warped product with the Lee-vector field $\lambda$ is orthogonal to $\mathfrak{D}^\perp$ if and only if
\begin{align}\label{3.1}
A_{JZ} X=g(J\lambda,JX)-(g(J\lambda,X)+JX(\mu))Z.
\end{align}
for some smooth function $\mu$ on $M$ satisfying $W\mu=0$ for all
$X\in\mathfrak{D},~~Z,W\in \mathfrak{D}^\perp$.
\section{Chen's second inequality for CR-warped produts}
In this section, we establish a sharp inequality for the second fundamental form of CR-warped products in LCK-space forms.

Let $M$ be a (pseudo-)Riemannian $k$-manifold with inner
product $g$ and $e_1,\cdots,e_k$ be an orthonormal frame fields on $M$.
For a differentiable function $\phi$ on $M$, the gradient $\nabla
\phi$ and the Laplacian $\triangle \phi$ of $\phi$ are defined
respectively by
\begin{align}\label{3.2}
&g(\nabla \phi,X)=X(\phi),\notag\\
&\triangle \phi=\sum\limits_{j=1}\limits^k\{(\nabla_{e_j} e_j)\phi-e_j e_j(\phi) \}=-div \nabla \phi
\end{align}
for vector field $X$ tangent to $M$, where $\nabla$ is the
Riemannian connection on $M$. As a consequence, we have
\begin{align}\label{3.3}
\|\nabla \phi\|^2=\sum\limits_{j=1}\limits^k (e_j(\phi))^2.
\end{align}
Using the above results, we will prove our main theorem.

\begin{theorem}\label{T:main} Let $M^n=N^T\times_f N^\perp$ be a $CR$-warped
product submanifold of an $LCK$-space form $\tilde{M}^{2m} (c)$. Then 
\begin{enumerate}
\item[(i)] The second fundamental form $h$ of $M$ satisfies the following inequality
\begin{align}\label{3.4}
\|h)\|^2\geq\frac{c}{2}-2p(\triangle \ln f)+4p\|\nabla \ln f\|^2-4pk (\lambda \ln f)-\frac{p+2k}{2}\|\alpha\|^2 
+\frac{1}{2} G + 2p F.
\end{align}
where 
\begin{align*}
G&=\sum\limits_{j=1}\limits^{k}\{g(\alpha^\#,e_j)^2+g(\alpha^\#,\nabla_{e_j}e_j)+g(\nabla_{e_j}\lambda ,e_j)-g(\alpha^\#,e_j)\nabla_{e_j}e_j\}~~~~~~~~~~~~~~~~\\
F&=\sum\limits_{j=1}\limits^k\{ \ (g(\beta^\#,e_j))^2+(g(\alpha^\#,e_j))^2 -g(\nabla_{e_j}\lambda ,e_j)-g(\alpha^\#,\nabla_{e_j}e_j)\notag\\
&-g(\nabla_{Je_j}\lambda ,Je_j)-g(\alpha^\#,\nabla_{Je_j}Je_j)-g(\beta^\#,\nabla_{Je_j} e_j)+g(\beta^\#,\nabla_{e_j} Je_j)\}
\end{align*}
and $\dim(N^T)=2k,\,\dim(N^\perp)=p$ and $\lambda$ is the Lee vector field orthogonal to $\mathfrak{D}^\perp$ in $M$.
\item[(ii)] If the equality sign holds in \eqref{3.4}, then $N^T$ is a totally geodesic submanifold of $\tilde M (c)$ and $N^\perp$ is totally umbilical in $\tilde M (c)$ with $h(\mathfrak{D}, \mathfrak{D})\perp\nu$ and $h(\mathfrak{D}^\perp, \mathfrak{D}^\perp)\perp J\mathfrak{D}^\perp$.
\end{enumerate}
\end{theorem}
\begin{proof} 
Let $M=N^T\times_f N^\perp$ be a $CR$-warped
product submanifold of an $LCK$-space form $\tilde{M}^{2m} (c)$. Then from the definition of $h$, we have
\begin{align}\label{S1}
\|h\|^2&=\|h(\mathfrak{D},\mathfrak{D})\|^2+2\|h(\mathfrak{D},\mathfrak{D}^\perp)\|^2+\|h(\mathfrak{D}^\perp,\mathfrak{D}^\perp)\|^2\notag\\
&=\|h_{J\mathfrak{D}^\perp}(\mathfrak{D},\mathfrak{D})\|^2+\|h_{\nu}(\mathfrak{D},\mathfrak{D})\|^2+2\|h(\mathfrak{D},\mathfrak{D}^\perp)\|^2\notag\\
&\hspace{.3cm}+\|h_{J\mathfrak{D}^\perp}(\mathfrak{D}^\perp,\mathfrak{D}^\perp)\|^2+\|h_{\nu}(\mathfrak{D}^\perp,\mathfrak{D}^\perp)\|^2.
\end{align}
The first term in the right hand side of \eqref{S1} vanishes identically by using Lemma \ref{L4} (i). Then, leaving the positive second, forth and fifth terms in \eqref{S1}, we get
\begin{align}\label{S2}
\|h\|^2\geq2\|h(\mathfrak{D},\mathfrak{D}^\perp)\|^2
\end{align}
Now, we calculate right hand side as follows
\begin{align}\label{3.5}
\|h(\mathfrak{D},\mathfrak{D}^\perp)\|^2=\sum\limits_{j=1}\limits^{k}\sum\limits_{i=1}\limits^p \|h(X_j,Z_i)\|^2, 
\end{align}
where $X_j$ for $\{j=1,\cdots,2k\}$ and $Z_\alpha$ for
$\alpha=\{1,\cdots,p\}$ are orthonormal frames on $N^T$ and $N^\perp$,
respectively. On $N^T$ we will consider a local orthonormal
frame, namely $\{e_j,J e_j\}$, where $\{j=1,\cdots,k\}$. We have to evaluate $\|h(X,Z)\|^2$ with $X\in\mathfrak{D}$ and $Z\in
\mathfrak{D}^\perp$. The second fundamental form $h(X,Z)$ is normal to $M$ so,
it splits into two orthogonal components
\begin{align}\label{3.6}
h(X,Z)=h_{ J\mathfrak{D}^\perp}(X,Z)+h_\nu (X,Z), 
\end{align}
where $h_{J\mathfrak{D}^\perp}(X,Z) \in J\mathfrak{D}^\perp$ and $h_\nu (X,Z)\in \nu$. So
\begin{align}\label{3.7}
\|h(X,Z)\|^2=\|h_{J\mathfrak{D}^\perp}(X,Z)\|^2+\|h_\nu (X,Z)\|^2. 
\end{align}
Now, we have
\begin{align}\label{3.9}
 \| h_{J\mathfrak{D}^\perp}(\mathfrak{D},\mathfrak{D}^\perp)\|^2=\sum\limits_{j=1}\limits^{k}\sum\limits_{i=1}\limits^p\{\|h_{J\mathfrak{D}^\perp}(e_j,Z_i)\|^2+ \|h_{J\mathfrak{D}^\perp}(J e_j,Z_i)\|^2\}
 \end{align}
for any $X\in\mathfrak{D}$ and $Z\in \mathfrak{D}^\perp$. Using \eqref{3.3} and Lemma \ref{L4}, after the computations, we can conclude that
\begin{align}\label{3.10}
\|h_{J\mathfrak{D}^\perp}(\mathfrak{D},\mathfrak{D}^\perp)\|^2&=2p\|\nabla \ln f\|^2+p\sum\limits_{j=1}\limits^{k}\{g(J\lambda,e_j)^2+g(\lambda,e_j)^2\}\notag\\
&-2p\sum\limits_{j=1}\limits^{k}\{(Je_j \ln f)g(\lambda,Je_j)+(e_j\ln f)g(\lambda,e_j)\}.
\end{align}
Consider the tensor field $\tilde{H}_B$. As we already have seen
\begin{align}\label{3.12}
\tilde{H}_B (X,Z)=g((\tilde\nabla_{JX} )h(X,Z)-(\tilde\nabla_X h)(JX,Z),JZ),
\end{align}
for any $X\in\mathfrak{D}$ and $Z\in \mathfrak{D}^\perp$. Using the definition of $\tilde\nabla h$, we obtain
\begin{align}\label{3.13}
\tilde{H}_B (X,Z)&=g({\nabla^\perp}_{JX} h(X,Z)-h(\nabla_{JX} X,Z)-h(X,\nabla_{JX} Z),JZ)\notag\\
&-g({\nabla^\perp}_{X} h(JX,Z)-h(\nabla_{X} JX,Z)-h(JX,\nabla_{X} Z),JZ).
\end{align}
In order to solve easily, we separate each term as follows
\begin{align*}
&T_1=g({\nabla^\perp}_{JX} h(X,Z),JZ),\;\;\;T_2=-g(h(\nabla_{JX} X,Z),JZ),\notag\\
&T_3=-g(h(X,\nabla_{JX} Z),JZ),\,\;\;\;T_4=-g({\nabla^\perp}_{X} h(JX,Z),JZ)\notag\\
&T_5=g(h(\nabla_{X} JX,Z),JZ),\;\;\;T_6=g(h(JX,\nabla_{X} Z),JZ).
\end{align*}
First we will compute $T_1$ and $T_4$
\begin{align*}
T_4&=-{Xg(h(JX,Z),JZ)+g(h(JX,Z),{\nabla^\perp}_{X}JZ)}\notag\\
&=-X[g((\tilde\nabla_{Z}J)X, JZ)+g(\tilde\nabla_{Z}X,Z)]+g(h(JX,Z),{J\tilde\nabla}_{X}Z+(\tilde\nabla_{X}J)Z).
\end{align*}
Then, using \eqref{2.2}, we get
\begin{align*}
T_4&=X(g(\lambda ,X)g(Z,Z))- X(X\ln f g(Z,Z))+g(Jh(X,Z),h(JX,Z))\notag\\
&+X\ln f g(h(JX,Z),JZ).
\end{align*}
Using Lemma \ref{L3}, after the computations, we can conclude
that 
\begin{align}
T_4&=\{Xlnf g(\lambda ,X)- X^2lnf -(Xlnf)^2
 +g(\nabla_{X}\lambda ,X)\notag\\
&+g(\lambda,\nabla_{X}X)\}\|Z\|^2 +g(Jh(X,Z),h(JX,Z)).
\end{align}
Similarly, we find
\begin{align}\label{3.16}
T_1&=\{JXlnfg(\lambda ,JX)-(JX)^2(ln f)-(JX\ln f)^2+g(\nabla_{JX}\lambda ,JX)\notag\\
&+g(\lambda,\nabla_{JX}JX)\}\|Z\|^2+g(Jh(X,Z),h(JX,Z)).
\end{align}
Then, it is not difficult to obtain the following terms 
\begin{align}\label{3.18}
T_2=\{(J\nabla_{JX}X)\ln f+g(J\lambda,\nabla_{JX} X) \}\|Z\|^2
\end{align}
and
\begin{align}\label{3.19}
T_5=-\{g(J\lambda,\nabla_{X} JX)+(J\nabla_{X}JX)\ln f\}\|Z\|^2.
\end{align}
We direct our attention to the third and sixth terms:
\begin{align}\label{3.20}
T_3=\{(JX\ln f)^2+g(J\lambda, X)(JX\ln f)\}\|Z\|^2
\end{align}
and
\begin{align}\label{3.21}
T_6=\{(X\ln f)^2-g(\lambda,X)(X\ln f)\}\|Z\|^2.
\end{align}
After using all above expressions, equation \eqref{3.13} becomes 
\begin{align}\label{3.22}
 \tilde{H}_B (X, Z)&=\|Z\|^2\{JX\ln f g(\lambda, JX)-(JX)^2(\ln f)-(JX\ln f)^2+g(\nabla_{JX}\lambda, JX)\notag\\
&+g(\lambda,\nabla_{JX}JX)+X \ln f g(\lambda, X)- X^2\ln f -(X\ln f)^2+g(\nabla_{X}\lambda, X)\notag\\
&+g(\lambda,\nabla_{X}X)+(J\nabla_{JX}X)\ln f+g(J\lambda,\nabla_{JX} X)-g(J\lambda,\nabla_{X} JX)\notag\\
&-(J\nabla_{X}JX)\ln f+(JX\ln f)^2+g(J\lambda, X)(JX\ln f)\notag\\
&+(X\ln f)^2-g(\lambda,X)(X\ln f)\}+2g(Jh(X, Z), h(JX,Z)).
\end{align}
On the other hand, we easily find the following relation
\begin{align}\label{3.23}
(J\nabla_{JX}X)\ln f&=(\nabla_{JX}JX)\ln f-(JX\ln f)g(\lambda,JX)-(X\ln f)g(\lambda ,X)\notag\\
&+(\lambda\ln f)\|X\|^2,
\end{align}
Interchanging $X$ by $JX$ in \eqref{3.23}, we get
\begin{align}\label{3.24}
(J\nabla_X JX)\ln f&=-(\nabla_X X)\ln f+(JX\ln f)g(\lambda, JX)+(X\ln f)g(\lambda, X)\notag\\
&-(\lambda\ln f)\|X\|^2.
\end{align}
Using \eqref{3.23}, \eqref{3.24} in \eqref{3.22}, we get
\begin{align}\label{3.25}
\tilde{H}_B (X, Z)&=\{(\nabla_{JX}JX-(JX)^2)\ln f+(\nabla_X X-(X)^2)\ln f -2(X ln f)g(\lambda,X)\notag\\
&-2(JX\ln f)g(\lambda,JX)+2\lambda \ln f |X\|^2
+g(\nabla_{X}\lambda ,X)
+g(\lambda,\nabla_{X}X)\notag\\
&+g(\nabla_{JX}\lambda ,JX)+g(\lambda,\nabla_{JX}JX)+g(J\lambda,\nabla_{JX} X)\notag\\
&-g(J\lambda,\nabla_{X} JX)\} \|Z\|^2+2g(Jh(X,Z),h(JX,Z)).
\end{align}
Using orthonormal frame fields and Lemma \ref{L5}, we obtain
\begin{align}\label{3.26}
\tilde{H}_B (e_j,Z_i)&=\{((\nabla_{Je_j}Je_j)-(Je_j)^2)\ln f+(\nabla_{e_j} e_j-(e_j)^2)\ln f\notag\\
&-2(e_j ln f)g(\lambda,e_j)-2(Je_j\ln f)g(\lambda,Je_j)+2\lambda \ln f |e_j\|^2\notag\\
&+g(\nabla_{e_j}\lambda ,e_j)+g(\lambda,\nabla_{e_j}e_j)+g(\nabla_{Je_j}\lambda ,Je_j)+g(\lambda,\nabla_{Je_j}Je_j)\notag\\
&+g(J\lambda,\nabla_{JX} e_j)-g(J\lambda,\nabla_{e_j} Je_j)\}\|Z_i\|^2+2\|h_{\nu} (e_j,Z_i)\|^2.
\end{align}
Similarly, we have
\begin{align}\label{3.27}
\tilde{H}_B (J e_j,Z_i)&=\{((\nabla_{e_j}e_j)-(e_j)^2)\ln f+(\nabla_{Je_j} Je_j-(Je_j)^2)\ln f\notag\\
&-2(J e_j ln f)g(\lambda,J e_j)-2( e_j\ln f)g(\lambda, e_j)+2\lambda \ln f | e_j\|^2\notag\\
&+g(\nabla_{J e_j}\lambda ,J e_j)
+g(\lambda,\nabla_{J e_j}J e_j)+g(\nabla_{ e_j}\lambda , e_j)+g(\lambda,\nabla_{ e_j} e_j)\notag\\
&-g(J\lambda,\nabla_{ e_j} J e_j)+g(J\lambda,\nabla_{J e_j}  e_j)\}\|Z_i\|^2+2\|h_{\nu} (Je_j,Z_i)\|^2.
\end{align}
On the other hand from \eqref{3.2}, we have
\begin{align}\label{3.28}
\triangle (\ln f)=&\sum\limits_{j=1}\limits^k\{(\nabla_{e_j} e_j)(\ln f)-{e_j}^2(\ln f) \}+\sum\limits_{j=1}\limits^k\{(\nabla_{Je_j} Je_j)(\ln f)-{Je_j}^2(\ln f) \}
\end{align}
and from \eqref{3.3}, we have
\begin{align}\label{3.29}
2\|\nabla \ln f\|^2=2\sum\limits_{j=1}\limits^k (e_j(\ln f))^2+2\sum\limits_{j=1}\limits^k (J e_j(\ln f))^2.
\end{align}
Taking the sum of \eqref{3.26} and \eqref{3.27} and using \eqref{3.28} and
\eqref{3.29}, we get
\begin{align}\label{3.30}
&2\sum\limits_{j=1}\limits^k\sum\limits_{i=1}\limits^p\{ |h_{\nu} (e_j,Z_i)\|^2+\|h_{\nu} (Je_j,Z_i)\|^2\}\\
&=-2p(\triangle \ln f)-4kp\lambda\ln f+\sum\limits_{j=1}\limits^k\sum\limits_{i=1}\limits^p\{\tilde{H}_B(e_j,Z_i)+\tilde{H}_B (Je_j,Z_i)\}
\notag\\
&+2p\sum\limits_{j=1}\limits^k\{2(e_j lnf)g(\lambda,e_j) +2(Je_j\ln f)g(\lambda,Je_j)-g(\nabla_{e_j}\lambda ,e_j)-g(\lambda,\nabla_{e_j}e_j)\notag\\
&-g(\nabla_{Je_j}\lambda ,Je_j)-g(\lambda,\nabla_{Je_j}Je_j)-g(J\lambda,\nabla_{Je_j} e_j)+g(J\lambda,\nabla_{e_j} Je_j)\}.\notag
\end{align}
Now, from eqref{2.8}, \eqref{3.7}, \eqref{3.10} and \eqref{3.30},  we derive
\begin{align}\label{3.31}
\|h(\mathfrak{D}, \mathfrak{D}^\perp)\|^2&= \frac{c}{4}-p(\triangle \ln f)+2p\|\nabla \ln f\|^2-2pk (\lambda \ln f)\notag\\
&-\frac{1}{4}\sum\limits_{j=1}\limits^{2k}\sum\limits_{i=1}\limits^p\{P(e_j,e_j)+P(Z_i,Z_i)\}\notag\\
&+p\sum\limits_{j=1}\limits^k\{ \ (g(J\lambda,e_j))^2+(g(\lambda,e_j))^2 -g(\nabla_{e_j}\lambda ,e_j)\notag\\
&-g(\lambda,\nabla_{e_j}e_j)-g(\nabla_{Je_j}\lambda ,Je_j)-g(\lambda,\nabla_{Je_j}Je_j)\notag\\
&-g(J\lambda,\nabla_{Je_j} e_j)+g(J\lambda,\nabla_{e_j} Je_j)\}.
\end{align}
Then, from \eqref{2.4}, the above equation takes the form
\begin{align}\label{3.32}
\|h(\mathfrak{D}, \mathfrak{D}^\perp)\|^2&= \frac{c}{4}-p(\triangle \ln f)+2p\|\nabla \ln f\|^2-2pk (\lambda \ln f)-\frac{p+2k}{4}\|\alpha\|^2 \notag\\
&+\frac{1}{4}\sum\limits_{j=1}\limits^{k}\{g(\lambda,e_j)^2+e_j(\lambda,e_j)-g(\lambda,e_j)\nabla_{e_j}e_j\}\notag\\
&+p\sum\limits_{j=1}\limits^k\{ \ (g(J\lambda,e_j))^2+(g(\lambda,e_j))^2 -g(\nabla_{e_j}\lambda ,e_j)\notag\\
&-g(\lambda,\nabla_{e_j}e_j)-g(\nabla_{Je_j}\lambda ,Je_j)-g(\lambda,\nabla_{Je_j}Je_j)\notag\\
&-g(J\lambda,\nabla_{Je_j} e_j)+g(J\lambda,\nabla_{e_j} Je_j)\}\notag\\
&= \frac{c}{4}-p(\triangle \ln f)+2p\|\nabla \ln f\|^2-2pk (\lambda \ln f)-\frac{p+2k}{4}\|\alpha\|^2 \notag\\
&+\frac{1}{4}\sum\limits_{j=1}\limits^{k}\{g(\lambda,e_j)^2+g(\lambda,\nabla_{e_j}e_j)+g(\nabla_{e_j}\lambda ,e_j)\notag\\
&-g(\lambda,e_j)\nabla_{e_j}e_j\}+p\sum\limits_{j=1}\limits^k\{ \ (g(J\lambda,e_j))^2+(g(\lambda,e_j))^2 -g(\nabla_{e_j}\lambda ,e_j)\notag\\
&-g(\lambda,\nabla_{e_j}e_j)-g(\nabla_{Je_j}\lambda ,Je_j)-g(\lambda,\nabla_{Je_j}Je_j)\notag\\
&-g(J\lambda,\nabla_{Je_j} e_j)+g(J\lambda,\nabla_{e_j} Je_j)\}.
\end{align}
Thus, from \eqref{S2} and \eqref{3.32}, we find 
\begin{align}\label{3.33}
\|h\|^2= \frac{c}{2}-2p(\triangle \ln f)+4p\|\nabla \ln f\|^2-4pk (\lambda \ln f)-\frac{p+2k}{2}\|\alpha\|^2 
+\frac{1}{2} G + 2p F.
\end{align}
where 
\begin{align*}
G&=\sum\limits_{j=1}\limits^{k}\{g(\lambda,e_j)^2+g(\lambda,\nabla_{e_j}e_j)+g(\nabla_{e_j}\lambda ,e_j)-g(\lambda,e_j)\nabla_{e_j}e_j\}\\
F&=\sum\limits_{j=1}\limits^k\{ \ (g(J\lambda,e_j))^2+(g(\lambda,e_j))^2 -g(\nabla_{e_j}\lambda ,e_j)-g(\lambda,\nabla_{e_j}e_j)\notag\\
&\hspace{.3cm}-g(\nabla_{Je_j}\lambda ,Je_j)-g(\lambda,\nabla_{Je_j}Je_j)-g(J\lambda,\nabla_{Je_j} e_j)+g(J\lambda,\nabla_{e_j} Je_j)\}.
\end{align*}
Hence, the inequality \eqref{3.4} follows from \eqref{3.33}. For the equality case, from the leaving and vanishing terms in \eqref{S1}, we find
\begin{align}\label{eq1}
h(\mathfrak{D}, \mathfrak{D})\perp J\mathfrak{D}^\perp,\;\;\;h(\mathfrak{D}^\perp, \mathfrak{D}^\perp)\perp \nu.
\end{align}
Also, from the hypothesis of theorem for equality case, we have
\begin{align}\label{eq2}
h(\mathfrak{D}, \mathfrak{D})\perp \nu,\;\;\;h(\mathfrak{D}^\perp, \mathfrak{D}^\perp)\perp J\mathfrak{D}^\perp.
\end{align}
Then, from \eqref{eq1} and \eqref{eq2}, we conclude that
\begin{align}\label{eq3}
h(\mathfrak{D}, \mathfrak{D})=0,\;\;\;h(\mathfrak{D}^\perp, \mathfrak{D}^\perp)=0.
\end{align}
Hence, $N^T$ is totally geodesic in $\tilde M$ due to the fact that $N^T$ is totally geodesic in $M$ \cite{Bi, Chen01}. Furthermore, from \eqref{eq3}, $N^\perp$ is totally umbilical in $\tilde M$ due to the fact that $N^\perp$ is totally umbilical in $M$ \cite{Bi, Chen01}.
\end{proof}
Now, using Theorem \ref{T:main}, we have the following result for Vaisman manifolds.
\begin{theorem}\label{C2}
Let $M^n=N^T\times_f N^\perp$ be a $CR$-warped product submanifold of a Vaisman
space form $\tilde{M}^{2m} (c)$. Then the second fundamental form of $M$
satisfies the following inequality
\begin{align}\label{S3}
\|h\|^2\geq \frac{c}{2}-2p(\triangle \ln f)+4p\|\nabla \ln f\|^2-4pk (\lambda \ln f)-\frac{p+2k}{2}\|\alpha\|^2+2p G^{*}.
\end{align}
where  $\alpha$ is nonzero constant 1-form, $\|\alpha\|^2$ is the length of Lee form with respect to $g$ and
\begin{align*}
G^{*}&=\sum\limits_{j=1}\limits^k\{ \ (g(\beta^{\#},e_j))^2+\frac{1+2p}{2p}(g(\alpha^{\#},e_j))^2 -g(\nabla_{e_j}\lambda ,e_j)\notag\\
&-g(\alpha^{\#},\nabla_{e_j}e_j)-g(\nabla_{Je_j}\lambda ,Je_j)-g(\alpha^{\#},\nabla_{Je_j}Je_j)\notag\\
&-g(\beta^{\#},\nabla_{Je_j} e_j)+g(\beta^{\#},\nabla_{e_j} Je_j)\}.
\end{align*}
\end{theorem}
\begin{proof} The proof follows from Theorem \ref{T:main}, \eqref{2.3} and the characteristic of Vaisman manifold.
\end{proof}

\begin{acknowledgements}
This project was funded by the Deanship of Scientific Research (DSR), King Abdulaziz University, Jeddah, Saudi Arabia under grant no. (KEP-PhD-88-130-38). The authors, therefore, acknowledge with thanks DSR technical and financial support.
\end{acknowledgements}

\end{document}